\documentclass[reqno]{amsart}
\usepackage{amssymb,graphics}
\usepackage[mathscr]{euscript}
\usepackage{epsfig}

\newcommand{\jsd}{join-sem\-i\-dis\-trib\-u\-tive}

\newtheorem{thm}{Theorem}[section] 
\newtheorem{cor}[thm]{Corollary}
\newtheorem{lm}[thm]{Lemma}

\newtheorem{prop}[thm]{Proposition}

\newtheorem{exm}[thm]{Example}
\newtheorem{rmk}[thm]{Remark}

\newtheorem{pbm}[thm]{Problem}

\theoremstyle{definition}
\newtheorem{df}[thm]{Definition}

\newcommand{\alg}[1]{{\mathbf #1}}        % an arbitrary algebra

\newcommand{\Co}{\alg{Co}(\mathbb{R}^n,X)}

\DeclareMathOperator{\Cl}{Cl}

\date{\today}

\begin{document}

\author{Kira Adaricheva}
\address{Harold Washington College, 30 East Lake St.,
Chicago, IL 60601, USA} \email{kadaricheva@ccc.edu}

\author{Marcel Wild}
\address{University of Stellenbosch, Private Bag XI, Matieland 7602,
South Africa} \email{mwild@sun.ac.za}

\title[Realization of abstract convex gemetries]{Realization of abstract convex geometries by point configurations. Part I.}

\begin{abstract}
The Edelman-Jamison problem is to characterize those abstract convex geometries that are representable by a set of points in the plane. We show that some natural modification of the Edelman-Jamison problem is equivalent to the well known $NP$-hard order type problem.
\end{abstract}

\maketitle
\footnotetext{While working on this paper, the first author was partially supported by INTAS grant N03-51-4110 ''Universal Algebra and Lattice Theory''}

\section{Introduction}
A finite closure space $(J, -)$ is called \emph{a convex geometry} (see, for example, \cite{EdJa}), if it satisfies \emph{the anti-exchange axiom}, i.e.
\[
\begin{aligned}
x\in\overline{A\cup\{y\}}\text{ and }x\notin A
\text{ imply that }y\notin\overline{A\cup\{x\}}\\
\text{ for all }x\neq y\text{ in }J\text{ and all closed }A\subseteq J.
\end{aligned}
\]
Given a closure space, one can associate with it the lattice of closed sets $\Cl (J,-)$; vice versa, every finite lattice $L$ represents the lattice of closed sets of a closure space defined on the set $J$ of join-irreducible elements of $L$. 
In particular, convex geometries correspond to locally lower distributive lattices which by definition are such that for each nonzero element $x$ the lattice generated by all lower covers of $x$ is Boolean. It is well known [9, p.19] that lower local distributivity is equivalent to the conjunction of lower semimodularity and join-semidistributivity. The latter property is defined by
$$(\forall x, y, z \in L) \quad (x\vee y = x \vee z\quad \Rightarrow\quad x \vee y = x \vee (y \wedge z))$$
Join semidistributivity is clearly inherited by sublattices, but lower distributivity generally isn't.

The following classical example of finite convex geometries shows how they earned their name.
Given a set of points $X$ in Euclidean space $\mathbb{R}^n$, one defines a closure operator on $X$ as follows: for any $Y \subseteq X$, $\overline{Y}= ch(Y) \cap X$, where $ch$ stands for \emph{the convex hull}. One easily verifies that such an operator satisfies the anti-exchange axiom. Thus, $(X,-)$ is a finite convex geometry. Denote by $\Co$ the closure lattice of this closure space, namely, the lattice of convex sets relative to $X$.

The current work was motivated by \cite{W1} and the following problem raised in \cite{AGT}: which lattices can be embedded
into $\Co$ for some $n \in \omega$ and some finite $X \subseteq \mathbb{R}^n$?
Is this the class of all finite \jsd\ lattices?
The positive answer for a proper subclass of \jsd\ lattice, namely, for all
finite lower bounded lattices, was given independently in \cite{Ad,WS}.

As one of the possible approaches to establish the structure of sublattices of $\Co$, one can ask about lattices {\it exactly} representable as $\Co$. Or, equivalently, what finite convex geometries can be realized as the convex sets relative to some finite point configurations in $n$-space?  

This is essentially the 
\begin{pbm}\label{EJP}
\bf{Edelman-Jamison Problem}\cite{EdJa}:
   
\emph{Characterize those convex geometries that are realizable by suitable point configuration in} $\mathbb{R}^n$.
\end{pbm}

We provide the partial solution to this problem for $n=2$ in the second part of our paper.
Here we discuss the connection of the Edelman-Jamison Problem to the $\bf{Order Type Problem}$. 

In combinatorial geometry, order types were introduced as a tool to capture essential features of point configurations. Assuming that a configuration is in a \emph{general position}, i.e. none of three distinct points are on one line, one defines an order-type of this configuration as a function of orientation of triples of distinct points. 

The Order Type Problem asks \emph{whether a given function from the triples of non-equal members from a given finite set $J$ into the two element set $\{-1,1\}$ can be realized as the orientation of triples of $|J|$ points on the plane in the general position}. It is known that the Order-Type problem is NP-hard. 

We show that point configurations that are equivalent as order types are also equivalent as convex geometries. On the other hand, there are plenty of point configurations that produce the same convex geometry while being non-equivalent as order types: see Example \ref{8points} and the follow-up series of examples described in Proposition \ref{series}. In fact, we show in Corollary \ref{plenty} that the number of non-equivalent order types corresponding to the same convex geometry cannot be polynomially bounded.

On the other hand, knowing the convex geometry formed by a given configuration, \emph{and the circular clock-wise order of the points} in the outside layer of this configuration, allows to determine the order type uniquely, see Theorem \ref{main}. Thus, for the convex geometries that enjoy a unique clock-wise circular order of their extreme points in each point configuration realizing them, the Edelman-Jamison Problem is polynomially equivalent to the Order-Type Problem.

\section{Convex 4-geometries}

Let $(X, \sim)$ be the convex geometry induced by a finite set of points $X \subseteq \mathbb{R}^2$ in general position. Thus its lattice of closed subspaces is $\alg{Co}(\mathbb{R}^2, X)$. A {\it rooted triangle of} $(X, \sim)$ is a pair $(T, \vec{x})$ such that $T \cup \{\vec{x}\} \subseteq X, \ |T| = 3, \ \vec{x} \in \widetilde{T} - T$ (thus $\vec{x}$ is in the interior of the triangle spanned by $T$).  Notice that a fixed $T$ may give rise to many rooted triangles $(T, \vec{x}), (T, \vec{y}), \cdots$ in $(X, \sim)$.
Since each polygon is partitioned by triangles, the closure operator $Y \mapsto \widetilde{Y}$ is determined\footnote{Put in other words, the family $\{T \rightarrow \{x\} | \ (T, x) \in \mathcal{R}\mathcal{T}(X)\}$ is an \emph{implicational base} in the sense of \cite{W2}.} by the set $\mathcal{R} \mathcal{T}(X)$ of all rooted triangles via
$$\widetilde{Y} \ = \ Y \cup \{\vec{x} \in X | \ \exists (T, \vec{x}) \in  \mathcal{R}\mathcal{T}(X) \ \mbox{{\it with}} \ T \subset Y\}.$$

A convex geometry $(J, -)$ which is isomorphic\footnote{Any two closure spaces are {\it isomorphic} if there is a bijection mapping one onto the other, while preserving the closure operator.} to the kind $(X, \sim)$ discussed above $(X \subseteq \mathbb{R}^2$ suitable), will be called  {\it realizable}. 

We need some preliminaries in order to formulate a necessary condition for being realizable. 

Call a subset $D$ of any closure space $(J, -)$  {\it 
dependent} if there is $x \in D$ with $x \in \overline{D - \{x\}}$. An inclusion-minimal dependent set $C$ is often called a {\it circuit} (adopting matroid terminology). It easily follows from the anti-exchange axiom that in a circuit $C$ of a convex geometry $(J, -)$ there is a {\it unique} element $x = x(C)$ of $C$, call it the {\it root} of $C$, such that $x \in \overline{C-\{x\}}$. 

It is not hard to show that one way to obtain a circuit with root $x$ is as follows. If $T \subseteq J -\{x\}$ is inclusion-minimal with $x \in \overline{T}$, then $C := T \cup \{x\}$ is a circuit. Let $Circ(J, -)$ be the set of all circuits of a convex geometry $(J, -)$. Thus, if $(J, -)$ happens to be realizable by some $X \subseteq \mathbb{R}^2$, then $Circ(J,-)$ bijectively corresponds to $\mathcal{R}\mathcal{T}(X)$ via $C \mapsto (C - \{x (C)\}, x (C))$.  

In particular, every realizable $(J, -)$ is a {\it convex 4-geometry} in that $|C| = 4$ for all $C \in Circ(J,-)$. It is handy to call a pair $(T, x), \ T$ any $3$-element set with $ x \not\in T$, a {\it quasi rooted triangle}. 

The following proposition is reminiscent of \cite{D} Theorem 7:

\begin{prop}\label{Dietrich} Let  $\mathcal{Q}\mathcal{R}\mathcal{T}$ be a family of candidate rooted triangles of a set $J$. Then the following are equivalent:
\begin{enumerate}
\item [(i)] There is a convex 4-geometry $(J, -)$ such that $(T, a) \mapsto T \cup \{a\}$ is a bijection from $\mathcal{Q}\mathcal{R}\mathcal{T}$ onto $Circ(J, -)$.
\vskip .5cm
\item[(ii)] Dietrich's axiom: For all $(T_1, a), (T_2, b) \in \mathcal{Q}\mathcal{R}\mathcal{T}$ with $a \in T_2$ there is  $(T_3, b) \in \mathcal{Q}\mathcal{R}\mathcal{T}$ with $T_3 \subseteq (T_1 \cup T_2) - \{a, b\}$. In words: Each triangle having a vertex ``colliding'' with another triangle's root, can be locally shifted whilst keeping its root (Figure 1).
\end{enumerate}

\begin{figure}[!h]
\begin{center}
\includegraphics[scale=.45]{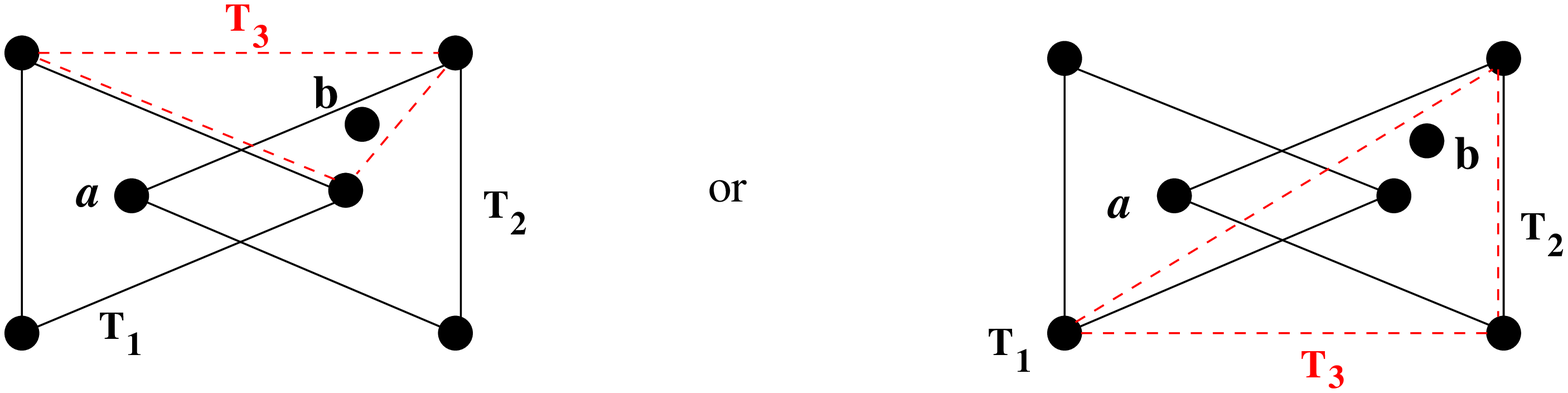}
\caption{}
\label{Dietrich1}
\end{center}
\end{figure}

\end{prop}

\begin{proof}: As to (i) $\Rightarrow$ (ii), let $T_1 \cup \{a\}$ and $T_2 \cup \{b\}$ be circuits of {\it any} convex geometry with roots $a, b$ respectively. Let $T:= (T_1 \cup T_2) - \{a, b\}$. Then $a \in \overline{T\cup \{b\}}$ and $b \in \overline{T \cup \{a\}}$. By the anti-exchange property either $a \in \overline{T}$ or $b \in \overline{T}$ takes place. In both cases $b \in \overline{T}$. Let $T_3 \subseteq T$ be minimial with $b \in \overline{T}_3$. Then, as mentioned previously, $T_3 \cup \{b\}$ is a circuit.

As to (ii) $\Rightarrow$ (i), we adhere to Figure 2 and first define
$$\overline{A}\  := \  A \cup \{a | \ \exists (T,a) \in \mathcal{Q}\mathcal{R}\mathcal{T} \ \mbox{with} \ T \subseteq A\}$$
Obviously this yields a monotone and extensive operator $\mathcal{P}(J) \mapsto \mathcal{P}(J)$. Suppose we had $\overline{\overline{A}} \neq \overline{A}$ for some $A \subseteq J$. Picking $b \in \overline{\overline{A}} - \overline{A}$ there would be some $(T_2, b)$ in $\mathcal{Q}\mathcal{R}\mathcal{T}$ with $T_2 \subseteq \overline{A}$ (and trivially $T_2 \not\subseteq A$). We may assume that among all possible $T_2$'s of this kind our $T_2$ minimizes $|T_2 - A|$.  Pick any $a \in T_2 - A$. Since $a \in \overline{A}$, there is a $(T_1, a)$ in $\mathcal{Q}\mathcal{R}\mathcal{T}$ with $T_1 \subseteq A$. By (ii) there is a $(T_3,b)$ in $\mathcal{Q}\mathcal{R}\mathcal{T}$ with $T_3 \subseteq (T_1 \cup T_2) - \{a, b\}$. Since $T_3 \subseteq \overline{A}$, and since $T_3 - A \subseteq T_2 - (A \cup \{a\})$ implies $|T_3-A| < |T_2-A|$, we get a contradiction to the minimality of $T_2$. Hence our operator $Y \mapsto \overline{Y}$ is idempotent, i.e. a closure operator. It is clear that the circuits with respect to this closure operator are precisely the members of $\mathcal{Q}\mathcal{R}\mathcal{T}$.

\begin{figure}[ht]
\begin{center}
\includegraphics[scale=.55]{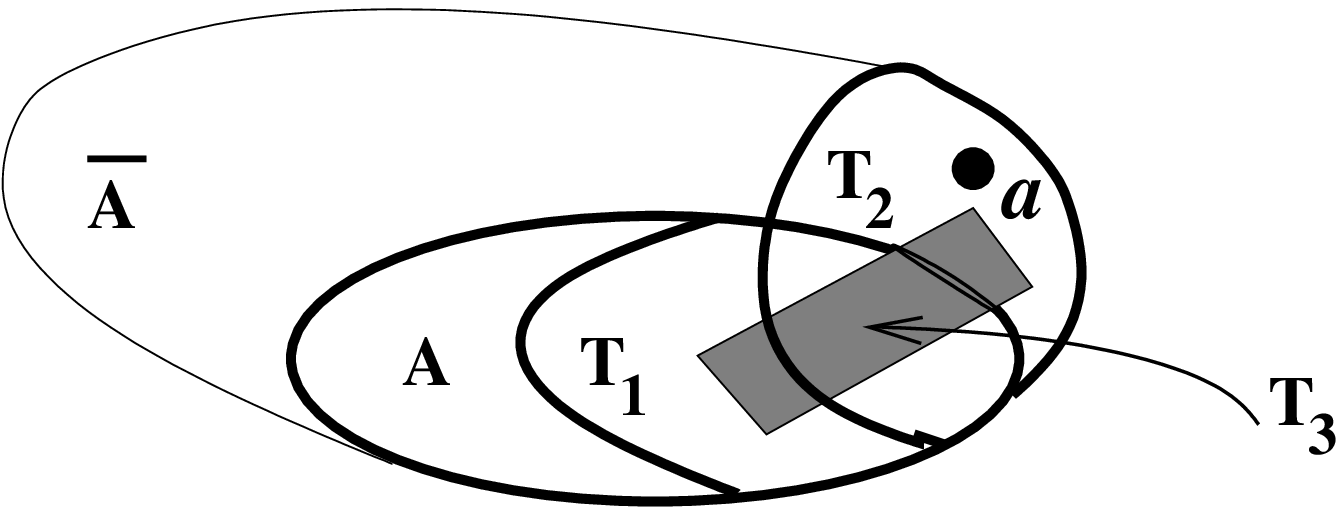}
\caption{}
\label{Dietrich2}
\end{center}
\end{figure}

In order to verify the anti-exchange property, suppose there was a closed $A \subseteq J$ and distinct elements $a, b \not\in A$ such that $a \in \overline{A \cup \{b\}}$ and $b \in \overline{A \cup \{a\}}$. Then there are $(T_1, a)$ and $(T_2, b)$ in $\mathcal{Q}\mathcal{R}\mathcal{T}$ with $T_1 \subseteq A\cup \{b\}$ and $T_2 \subseteq A\cup \{a\}$. By (ii) there is a  $(T_3, b)$ in $\mathcal{Q}\mathcal{R}\mathcal{T}$ with $T_3 \subseteq (T_1 \cup T_2) - \{a, b\}$. This implies $T_3 \subseteq A$, and whence the contradiction $b \in \overline{A} = A$. 
\end{proof}
In view of Proposition \ref{Dietrich} we adopt from now on the notation $(J,\mathcal{R}\mathcal{T})$ rather than $(J, -)$ for convex $4$-geometries. Here $\mathcal{R} \mathcal{T}$ is a set of rooted triangles (based on $J$) which satisfies Dietrich's axiom.

\subsection{Layers of convex $4$-geometries}
\vskip 0.5cm

For any finite convex geometry $(J,-)$ one can define, recursively, the family of subsets $L_i\subseteq J$, $i=0,1,\dots$, called \emph{layers}. Let $L_0= \{x \in J: x \not\in \overline{J\backslash\{x\}}\}$ be the set of extreme points of $(J,-)$.

Assume now that the layers $L_0,\dots,L_{n-1}$ are defined and $J\not = \bigcup_{i < n} L_i$.
Let $J_n=J \backslash \bigcup_{i < n} L_i$. There is a naturally defined convex geometry
$(J_n,-)$ whose closure operator is the restriction of closure operator of $(J,-)$ on $J_n$.
Then $L_n= \{x \in J_n: x \not\in \overline{J_n\backslash\{x\}}\}$ is the set of extreme points of geometry $(J_n,-)$.

One can proceed with defining the layers $L_0, \dots, L_k$ until $J \backslash \bigcup_{i\leq k} L_i = \emptyset$. This defines \emph{the complete family of layers} of $(J,-)$. We will call the number $q(J):=k$ \emph{the depth} of convex geometry. Besides, layer $L_0$ is called the \emph{outermost layer}, $L_k$ is the \emph{innermost layer}.

We collect the easy facts about layers in the following statement:

\begin{prop}
\begin{enumerate} Let $L_0,\dots, L_k$ be the complete family of layers of the convex geometry $(J,-)$.
\item[(1)] $L_i \cap L_j=\emptyset$, when $i\not = j$;
\item[(2)] $\bigcup_{i \leq k} L_i = J$;
\item[(3)] $L_{i+1}, \dots, L_k \subseteq \overline{L_i}$, for any $i < k$.
\end{enumerate}
\end{prop}
The proof follows easily from the well-known fact about convex geometries that $J=\overline{L_0}$.

\section{Order types}

For all noncolinear $\vec{x},\vec{y}, \vec{z} \in \mathbb{R}^2$ define
$$\mbox{sign} (\vec{x},\vec{y},\vec{z})\quad : =\quad \left\{ \begin{array}{rll} 1, & \mbox{if} & \vec{x},\vec{y},\vec{z} \ \mbox{are positively oriented (anticlockwise)}\\
-1, & \mbox{if} & \vec{x},\vec{y},\vec{z} \ \mbox{are negatively oriented} \end{array}\right.$$

Of course, that could be defined in terms of determinants, but there is no need to do so. Recall that  $\mbox{sign} (\vec{x}, \vec{y}, \vec{z}) = \ \mbox{sign} (\vec{z}, \vec{x}, \vec{y}) = \ \mbox{sign} (\vec{y}, \vec{z}, \vec{x})$ (cyclic permutability) and $\mbox{sign} (\vec{x}, \vec{y}, \vec{z}) = - \mbox{sign} (\vec{x}, \vec{z}, \vec{y})$.
\vspace{0.4cm}

For a set $J$ let $J[3]$ be the set of all triplets $(a,b,c)$ with distinct $a, b, c$ in $J$. Following \cite{AAK} we call two point configuration $X, Y \subseteq \mathbb{R}^2$ {\it equivalent} if there is a bijection $F: X \rightarrow Y$ which preserves the orientation of all triples in $X$. 
Following [7], call  $t: J[3] \rightarrow \{1, -1\}$ an {\it order type on} $J$,  
if there is a function $f: J \rightarrow \mathbb{R}^2$ such 
that for all $(a,b,c)$ in $J[3]$ one has
$$t(a,b,c) = \ \mbox{sign}(f(a), f(b), f(c))$$

The point configuration $X: = f(J)$ is then said to {\it realize} the order type $t$. In brief, $t$ is an order type, if it represents the orientation of triples of some suitable point configuration. If a particular $f$ is relevant, we shall write $t=t_f$.

Reminiscent to equivalent point configurations we declare two order types $t_1$ and $t_2$ on $J$ {\it equivalent}, if there is a bijection $\delta : \ J \rightarrow J$ such that $t_2 = t_1 \circ \overline{\delta}$. Here $\overline{\delta}: \ J [3] \rightarrow J[3]$ is the canonic map induced by $\delta$. A minute's thought confirms:
\begin{cor} Two order types $t_1$ and $t_2$ on $J$ are equivalent if and only if any two corresponding realizing point configurations $X_1$ and $X_2$ are equivalent.
\end{cor}
Let us call two order types $t_1, t_2$ {\it weakly equivalent} if $t_1$ is equivalent to either $t_2$ or $-t_2$. Accordingly weakly equivalent point configurations are defined. 
\begin{prop}\label{isomorphic}
Any two weakly equivalent point configurations $X,Y \subseteq \mathbb{R}^2$  induce isomorphic convex geometries.
\end{prop}
\begin{proof} Let $\vec{x}_1, \vec{x}_2, \vec{x}_3, \vec{x} \in X$ be distinct. A quick sketch confirms that
\vskip .4cm
(1) \quad $\vec{x} \in ch (\{\vec{x}_1, \vec{x}_2, \vec{x}_3\}) \quad \Leftrightarrow$\\
 
\hspace*{1.5cm} $\mbox{sign} (\vec{x}_1, \vec{x}_2, \vec{x}_3)  = \ \mbox{sign} (\vec{x}, \vec{x}_2, \vec{x}_3) = \ \mbox{sign} (\vec{x}_1, \vec{x}, \vec{x}_3)
 = \ \mbox{sign} (\vec{x}_1, \vec{x}_2, \vec{x}).$
\vskip .4cm
Thus, if $F: \ X \rightarrow Y$ is either orientation preserving or orientation reversing, then 
\vskip .6cm
$\vec{x} \in ch (\{\vec{x}_1, \vec{x}_2, \vec{x}_3\})\quad \Leftrightarrow$\\

$\mbox{sign} (\vec{x}_1, \vec{x}_2, \vec{x}_3) = \mbox{sign} (\vec{x}, \vec{x}_2, \vec{x}_3) = \mbox{sign} (\vec{x}_1, \vec{x}, \vec{x}_3) =
\mbox{sign} (\vec{x}_1, \vec{x}_2, \vec{x})  \quad \Leftrightarrow$\\

$ \mbox{sign}
(F\vec{x}_1, F\vec{x}_2, F\vec{x}_3) = \mbox{sign} (F \vec{x}, F\vec{x}_2, F\vec{x}_3) = \mbox{sign}(F\vec{x}_1, F\vec{x}, F\vec{x}_3)$$
 = \mbox{sign} (F \vec{x}_1, F\vec{x}_2, F\vec{x})$ \\

$\Leftrightarrow  F(\vec{x}) \in ch (\{F(\vec{x}_1), F(\vec{x}_2), F(\vec{x}_3)\}).$
\vskip .6cm
This shows that $F$ is an isomorphism of convex geometries.
\end{proof}

While point configuration are less abstract than order types, the latter will be more convenient in the proofs. Here is an appetizer.
\begin{prop} For each convex $4$-geometry  $(J, \mathcal{R}\mathcal{T})$  the following are equivalent:
\begin{enumerate}
 \item [(i)] $(J, \mathcal{R}\mathcal{T})$ is realizable
\item[(ii)] There is an order type $t$ on $J$ such that for all distinct $a, b, c, d$ in $J$ one has:
$$(\{a, b, c\}, d) \in \mathcal{R}\mathcal{T} \quad  \Leftrightarrow \quad t(a, b, c) = t(d, b, c) = t(a, d, c) = t(a, b, d)$$
\end{enumerate}
\end{prop}
\begin{proof} 

(i) $\Rightarrow$ (ii). If $f: \ J \rightarrow \mathbb{R}^2$ is a realization of $(J,\mathcal{R}\mathcal{T})$, then $t: = t_f$ does the job since for all $a, b, c, d \in J$:
\vskip .3cm
\hspace*{1cm} $(\{a, b, c\}, d)\in  \mathcal{R}\mathcal{T}\quad  \Leftrightarrow \quad f(d) \in ch (\{f(a), f(b), f(c)\})\quad \stackrel{(1)}{\Leftrightarrow}$
\vskip .3cm
\hspace*{4cm} $t(a,b,c) = t(d,b,c) = t(a,d,c) = t(a,b,d)$
\vskip .3cm
(ii) $\Rightarrow$ (i). Let $t = t_f$ be an order type as in (ii). Then $f: \ J \rightarrow \mathbb{R}^2$ is a realization of $(J, \mathcal{R}\mathcal{T})$ because for all distinct $a, b, c, d$ in $J$ one has:

\hspace*{1cm} $(\{a, b, c\}, d) \in \mathcal{R}\mathcal{T} \ \Leftrightarrow \ t(a,b,c) = t(d,b,c) =t(a,d,c) =t(a,b,d)\quad \stackrel{(1)}{\Leftrightarrow}$
\vskip .3cm
\hspace*{4cm} $f(d) \in ch (\{f(a), f(b), f(c)\})$. 
\end{proof}

Let us illustrate these concepts with two examples.

\begin{exm} Consider these two point configurations:
\end{exm}
\begin{figure}[ht]
\begin{center}
\includegraphics[scale=.45]{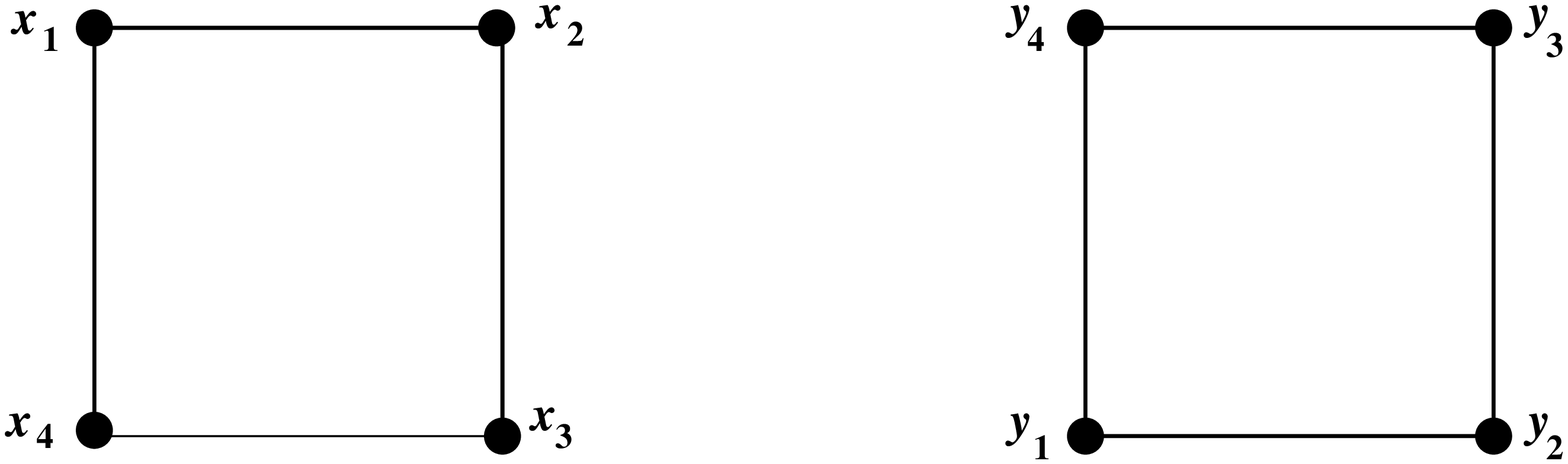}
\caption{}
\label{squares1}
\end{center}
\end{figure}

The labeling suggests that square $Y$ arises from square $X$ by reflecting the latter on its horizontal axis of symmetry. Thus, since $F(x_i) := y_i$ is an order reversing bijection, the point configurations $X$ and $Y$ are weakly equivalent. Rephrasing it in terms of order types, put $J := \{1, 2, 3, 4\}$ and define $t_i : J[3]  \rightarrow \{-1,1\} \ (i = 1, 2)$ by
$$t_1(i, j,k) := \mbox{sign} (x_i, x_j, x_k), \quad t_2 (i, j, k):= \mbox{sign} (y_i, y_j, y_k)$$
Then $t_2 = -t_1$, i.e. $t_1, t_2$ are weakly equivalent order types (putting $\delta = id$). Of course $X$ and $Y$ (as well as the corresponding order types) are actually equivalent since substituting $F$ by the bijection
$$G(x_1): = y_4, \quad G(x_2):= y_3, \quad G(x_4):= y_1, \quad G(x_3):=y_2$$
does the job.

Let us now fill the squares with six points each:
\begin{figure}[ht]
\begin{center}
\includegraphics[scale=.45]{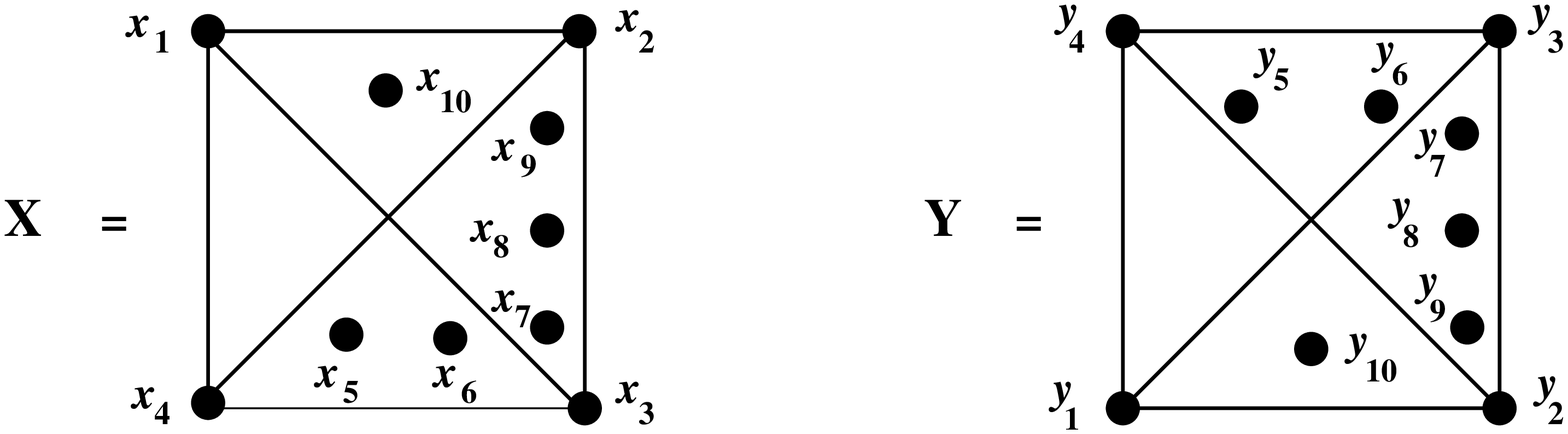}
\caption{}
\label{squares2}
\end{center}
\end{figure}

Again $Y$ is obtained from $X$ upon reflection on the middle axis. Put another way, looking at $X$ from ``below the sheet'' yields $Y$. Hence the two point configurations are again weakly equivalent. However, this time we will not succeed in finding a $G: X \rightarrow Y$ that establishes the equivalence of $X$ and $Y$. It suffices to show that 
\vskip .4cm
(2) \hspace*{2cm} $G(x_i) = y_i \ \mbox{for any equivalence}  \ G : X \rightarrow Y$
\vskip .6cm
because then e.g. $\mbox{sign} (y_1, y_2, y_3) =  1 \neq -1 = \mbox{sign} (x_1, x_2, x_3)$. In order to see (2), recall from Proposition \ref{isomorphic} that an equivalence $G: X \rightarrow Y$ is an isomorphism of convex $4$-geometries, i.e. preserves rooted triangles. Because the number of roots inside the triangles $\{x_1, x_2, x_3\}, \{x_1, x_2, x_4\}, \{x_1, x_3, x_4\}, \{x_2, x_3, x_4\}$ is $4, 1, 2, 5$ respectively (and dito for $Y$), it follows that $G$ maps $\{x_1, x_2, x)$ onto $\{y_1, y_2, y_3\}$, and so on. A moment's thought confirms that this forces (2).

The following example shows that there exist point configurations that are not weakly equivalent yet yield isomorphic convex $4$-geometries.

\newpage
\begin{figure}[ht]
\begin{center}
\includegraphics[scale=.4]{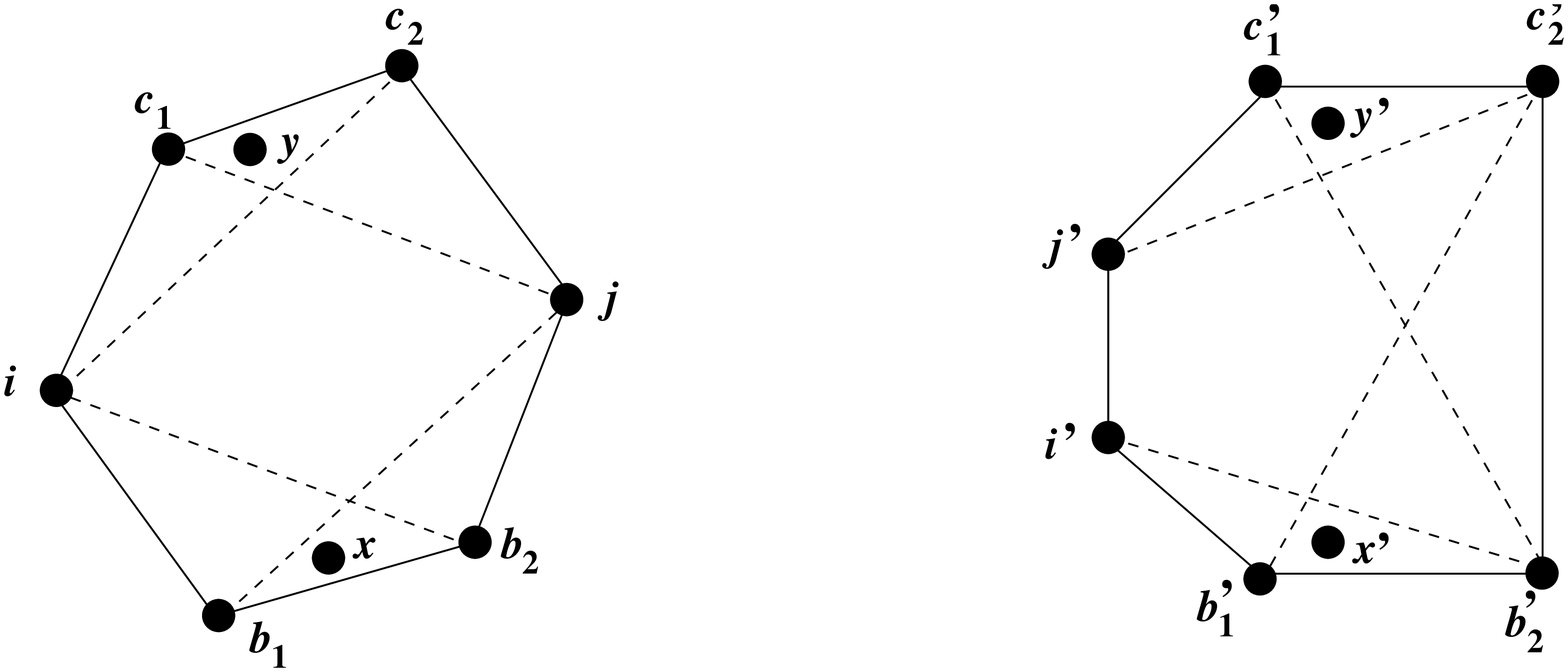}
\caption{}
\label{not weak}
\end{center}
\end{figure}

\begin{exm}\label{8points}
{}
\end{exm}
Let $L=\{b_1,b_2,c_1,c_2,i,j,x,y\}$ and $L'=\{b_1',b_2',c_1',c_2',i',j',x',y'\}$ be two $8$-point configurations given on Fig.1. It is easy to verify that the correspondence $s \rightarrow s'$  induces an isomorphism of the convex $4$-geometries defined  by $L$ and $L'$. Indeed, there are exactly $10$ rooted triangles in $L$ (correspondingly, in $L'$ after replacement of each $s$ by $s'$) : 
\vskip .6cm
$(\{ i,c_1,c_2\},y)$,          $(\{ i,b_1,b_2\},x)$,
$(\{b_1,c_1,c_2\},y )$,        $(\{ c_1,b_1,b_2\},x )$,
$(\{ b_2,c_1,c_2\},y  )$,  
\vskip .4cm
$(\{ c_2,b_1,b_2\},x )$,
$(\{ j,c_1,c_2\},y )$,          $(\{ j,b_1,b_2\},x )$,
$(\{y, b_1, b_2\}, x )$, $(\{c_1, c_2, x\}, y)$.
\vskip .6cm
On the other hand, these two point configurations are not weakly equivalent. Indeed, consider the following property of an extreme point $s$ 
in configuration $L$:

Among all the lines through $s$ and the other extreme points, exactly three separate $x$ and $y$.

There are only two points $s$ among $\{
b_1,b_2,c_1,c_2,i,j\}$ with this property, namely,  $i$ and $j$. 

In configuration $L'$, there are also only two points among $\{
b_1',b_2',c_1',c_2',i',j'\}$ with the property of separating $x'$ and $y'$, namely $c_2'$ and $b_2'$.

Crucially, the property that a line, say {\it line}$(c_2', i')$, separates $x'$ and $y'$, can be expressed in terms of orientations: $\mbox{sign} (c'_2, i', x') \neq \mbox{sign}(c'_2, i', y')$. Hence, if $L$ and $L'$ would be weakly equivalent as order types, $i,j$ would need to be mapped to $c_2'$ and $b_2'$. Also, being inner points, $x,y$ would need to be mapped to $x',y'$.
But such a mapping cannot preserve, neither reverse, the order type since {\it line}$(i,j)$ separates $x,y$ and {\it line}$(c_2',b_2')$ does not separate $x',y'$.

\section{Order types of a realizable convex $4$-geometry}

This leads us to define {\it Order-Types}$(J, \mathcal{R}\mathcal{T})$ as the set of all order types $t_f: \ J[3] \rightarrow \{-1,1\}$ induced by realizations $f: \ J \rightarrow \mathbb{R}^2$ of the convex $4$-geometry $(J, \mathcal{R}\mathcal{T})$. In particular, {\it Order-Types}$(J,\mathcal{R}\mathcal{T}) = \emptyset$ if $(J,\mathcal{R}\mathcal{T})$ is not realizable. 
Let \underline{{\it Order-Types}}$(J, \mathcal{R}\mathcal{T})$ be the set of {\it equivalence classes} of order types contained in {\it Order-Types}$(J, \mathcal{R}\mathcal{T})$.

Let us first dispense with the case of a {\it free} convex $4$-geometry  $(J, \emptyset)$. Observe that $q(J, \mathcal{R} \mathcal{T}) =1 \Leftrightarrow \mathcal{R}\mathcal{T} = \emptyset$. If  $|J| = n$ then any $n$-gon ($n$ points in general position) is a realization of $(J, \emptyset)$. Furthermore, if $X, Y \subseteq \mathbb{R}^2$ are two $n$-gons, then there are $n$ bijections $F: \ X \rightarrow Y$ that preserve the orientations of triples, namely precisely those $F$'s that map a fixed cyclic (say clockwise) ordering of $X$ onto one of the $n$ cyclic clockwise enumerations of $Y$. In particular\footnote{As an exercise, prove that $|\mbox{Order-Types} (J, \emptyset)|= (n-1)!$. In the present article (i.e. Part I) we stick to \underline{Order-Types}. More about Order-Types and automorphisms of 4-geometries follows in Part II.}, 
$|\underline{Order\mbox{-}Types}(J,\emptyset)| =1$.

\vskip 1cm
\subsection{Realizable convex $4$-geometries with few order types}

Consider the case of a convex $4$-geometry $(J,\mathcal{R}\mathcal{T})$ with $q(J)=2$ with just {\it one} interior point. Denote by $L$ the outside layer of $(J,\mathcal{R}\mathcal{T})$, and let $p$ be the unique point of the inside layer. 

Underlying the results of Edelman and Larman \cite{EdLa}, is the definition of \emph{equivalent elements} of $L$: Put $s \equiv t$, if 
$$(\forall u, v \in L) \ (\{u, v, s\}, p) \in \mathcal{R}\mathcal{T} \Leftrightarrow (\{u, v, t\}, p) \in \mathcal{R}\mathcal{T}.$$

In particular, $s,t$ cannot be equivalent, when $s,t$ are in a common rooted triangle.

It turns out that if $s_1,s_2, \dots, s_k$ are equivalent elements in some realizable convex $4$-geometry $(J,\mathcal{R}\mathcal{T})$, then in any realization of $(J,\mathcal{R}\mathcal{T})$, all elements $s_1,s_2, \dots, s_k$ appear in one \emph{cluster} in a circular order of layer $L$, i.e. no proper subset of $\{s_1,s_2, \dots, s_k\}$ can be flanked in that circular order by some elements non-equivalent to $s_i$.

This observation, even though not spelled out in \cite{EdLa}, led to the notion of \emph{irreducible} convex $4$-geometry $(J,\mathcal{R}\mathcal{T})$ as one with no equivalent elements in $L$. The idea was to reduce clusters of equivalent elements of the outside layer $L$ to unique points and consider this simplified convex $4$-geometry about which one can make strong statements. One of the crucial results of \cite{EdLa} is 
\begin{thm}\label{one-point}(Theorem 3.5 in \cite{EdLa})
If $(J,\mathcal{R}\mathcal{T})$ is a realizable irreducible convex $4$-geometry with the outside layer $L$ and one inner point then the circular order of $L$ in any point realization is determined uniquely up to reflection.
\end{thm}

\begin{cor} If $(J,\mathcal{R}\mathcal{T})$ is a realizable convex $4$-geometry with outside layer $L$ and one inner point, then {\it \underline{Order-Types}}$(J, \mathcal{R}\mathcal{T})$ has cardinality at most two. 
\end{cor}
\begin{proof} According to Theorem \ref{one-point} above, there is maximum two non-equivalent order types of irreducible convex geometry $(J',\mathcal{R}\mathcal{T})$  deduced from $(J,\mathcal{R}\mathcal{T})$. Since all equivalent elements of $L$ are located in clusters, and since any order of equivalent elements within a cluster produces an equivalent order type, we can obtain at most two nonequivalent order types $t_1, t_2$. In fact, for self-symmetric point configurations, $t_1$ and $t_2$ may be equivalent order-types. 
\end{proof}

We can  follow-up with the definition of equivalent elements of the outermost layer $L_0$ of a convex $4$-geometry $(J, \mathcal{R}\mathcal{T})$ in general. If $p$ is a point inside that layer, then we can define an equivalence $\equiv_p$ on $L_0$ as follows:
$s\equiv_p t$, if $(\{s,u,v\},p)$ is a rooted triangle iff $(\{t,u,v\},p)$ is a rooted triangle, for any $u,v \in L_0$. If $P$ is the collection of all points $p$ inside layer $L_0$, then we can define $s\equiv t$ and call $s,t$ \emph{equivalent} iff $s \equiv_p t$, for every $p \in P$. 

It turns out that, unlike the case of one inner point, the equivalent elements of a layer no longer should appear in clusters, even when we increase the number of inner points by just one. Return to Example \ref{8points} for an illustration. One can directly check that points $i,j$ are equivalent points on the outside layer of the convex $4$-geometry, and they appear in one cluster in realization $L'$, while they are flanked on both sides by non-equivalent points in point configuration $L$.

\begin{df}\label{simple}
We will call $(J,\mathcal{R}\mathcal{T})$ \emph{simple}, if the outermost layer $L_0$ of $(J,\mathcal{R}\mathcal{T})$ does not have  equivalent elements. 
\end{df}

Evidently, the notion of a simple geometry corresponds to ''irreducible'' convex geometry of \cite{EdLa} in case of one inner point.
 
\begin{thm}\label{ordering} 
In any realization $f: \ J \rightarrow \mathbb{R}^2$ of a realizable simple convex $4$-geometry $(J, \mathcal{R}\mathcal{T})$ with $q(J,\mathcal{R}\mathcal{T}) \geq 2$, the cyclic ordering of the outermost layer is uniquely determined up to reflection. 

\end{thm}
\begin{proof} Let $L_0$ be any outside layer of  $(J,\mathcal{R}\mathcal{T})$. Let $P$ be the set of inner points for this layer. For any $p \in P$, consider a sub-geometry of $(J,\mathcal{R}\mathcal{T})$ defined on $J_p=L \cup \{p\}$. According to Theorem \ref{one-point}, there exists unique up to reflection circular order of clusters $S_1,\dots, S_{k_p}$ of $\equiv_p$-equivalent elements of layer $L_0$. We claim that this ''partial'' circular order of $L$ can be uniquely extended to a ``linear'' circular order of $L_0$.

If there exists a cluster of more than one point, say, $S_1$ has points $s_1,s_2$, then there should be an inner point $q$ such that $(s_1,s_2) \not \in \equiv_q$. Since $(J,\mathcal{R}\mathcal{T})$ is realizable, there should be a circular order of clusters of $\equiv_q$-equivalent points compatible with $S_1,\dots,S_k$ (we choose one of two existing that will follow the orientation of the first choice). The intersection of these two equivalences on $L_0$ will provide a new equivalence, i.e. a linear order of clusters, in which $s_1,s_2$ are no longer in one cluster.
\end{proof}

\subsection{Realizable convex $4$-geometries with many order types}

Example \ref{8points} gives an idea of series of examples of convex $4$-geometries with the growing number of non-equivalent order types. In the notation of Proposition below, $g(p)=\mathcal{O}(p^k)$ will mean that $0<\lim_{p \to \infty}\frac{g(p)}{p^k}<\infty$.
\newpage
\begin{figure}[ht]
\begin{center}
\includegraphics[scale=.45]{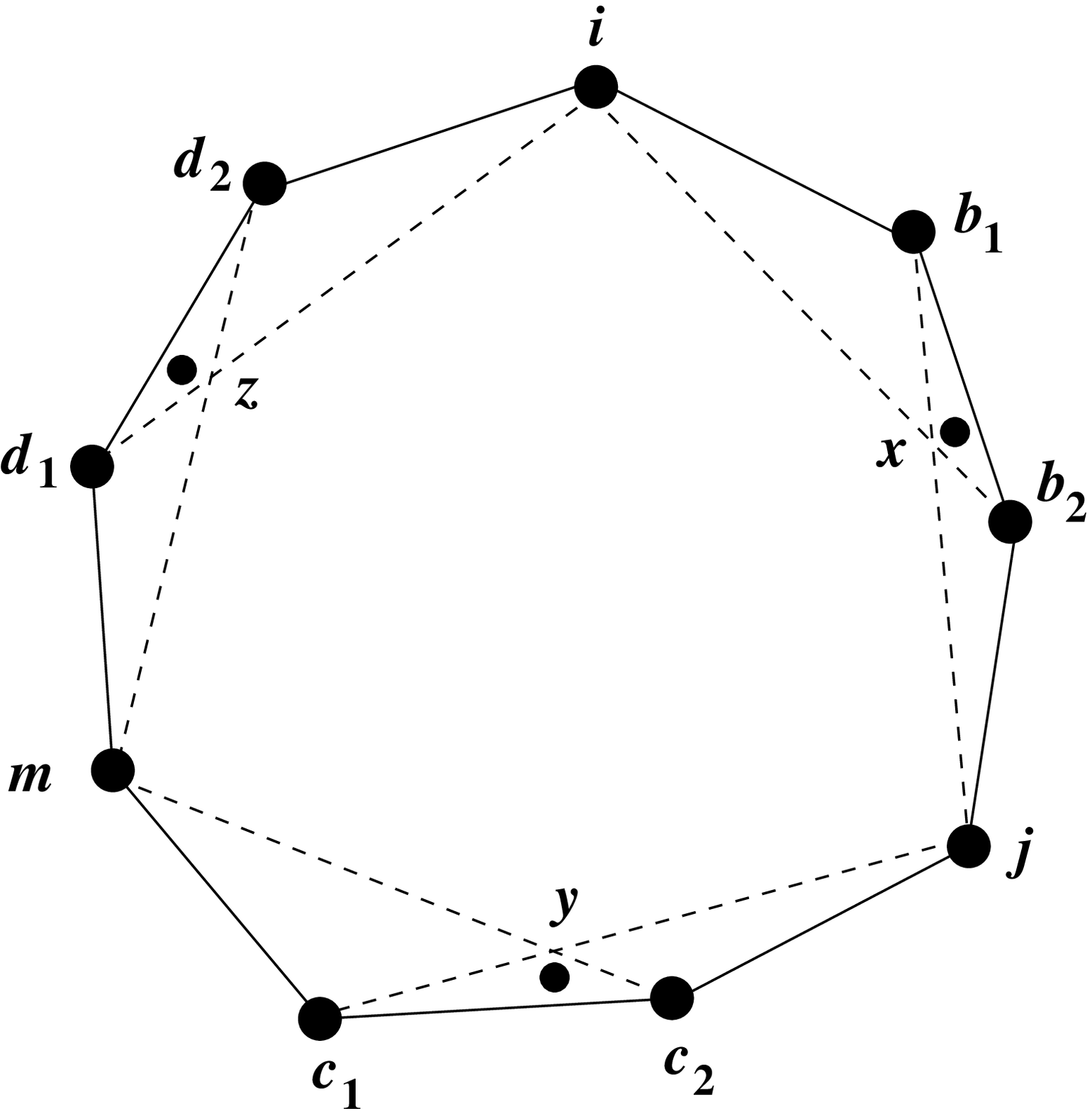}
\caption{}
\label{triangle}
\end{center}
\end{figure}

\begin{prop}\label{series} For any $k \in \omega$, there exists a series $\mathcal{J}_k=\{J(p): p > k \}$ of realizable convex $4$-geometries $J(p)$ with two layers such that
$f(p)=|J(p)|=\mathcal{O}(p)$, while $g(p)=|\text{{\it \underline{Order-Types}}}(J(p))|=\mathcal{O}(p^k)$.
\end{prop}
\begin{proof}
We first explain the idea for $k=1$ and $k=2$. When $k=1$, we use point configuration $L$ given on the left side of Figure \ref{not weak}. For any $p > 1$, let $J(p)$ be a convex $4$-geometry with two layers, whose first layer consists of $c_1,c_2,b_1,b_2$ and $p$ additional points, while the inner layer consists of two points $x$ and $y$. This convex geometry can be obtained from its point realization when $i,j$ of point configuration $L$ are replaced by $p'$ and $(p-p')$ points, correspondingly, for some $p' \leq p$. One needs to make sure that these points are just slightly displaced from original positions of $i$ and $j$ so that all $p$ points are in the same outside layer and they are equivalent.
When $p'=p$ one gets a configuration similar to $L'$ on the right side of Figure 
\ref{not weak}, where all $p$ points could be located on the segment between $i$ and $j$. Changing $p'$ from $p$ to $0$ one obtains $\left\lfloor \frac{p+1}{2}\right\rfloor$ different order types corresponding to the same convex geometry. Thus, $|\text{{\it \underline{Order-Types}}}(J(p))|=\mathcal{O}(p)$.

For $k=2$, consider the points configuration on Figure \ref{triangle}. It replicates point configuration $L$, replacing two equivalent elements of outside layer $i,j$ by three equivalent elements $i,j,m$ placed at the vertices of equilateral triangle. We then place $b_1,b_2$ on the arc connecting $i,j$, $c_1,c_2$ on the arc $j,m$ and $d_1,d_2$ on the arc $m,i$. The inner layer now consists of three elements $x,y,z$, placed close to the center of segments $[b_1,b_2]$, $[c_1,c_2]$, $[d_1,d_2]$, correspondingly. The number of rooted triangle is $21$: say, $x$ is inside $T(t,b_1,b_2)$, where $t$ ranges over all points of the outside layer other than $b_1,b_2$; similarly, for $y$ and $z$. For any $p > 2$ one could split $p$ points into three subsets of $p''$, $p'$ and $(p-p'-p'')$ elements, for some $p'+p'' \leq p$ and place, slightly displaced, into positions of $i,j$ and $m$ correspondingly, so that all $p$ points are equivalent. This produces the point realization for $J(p)$ from series $\mathcal{L}_3$. Evidently, $|J(p)|=\mathcal{O}(p)$. Varying $p''$ and $p'$ we may obtain about $\frac{1}{3}C_{p+1}^2$ of non-equivalent order-types, in particular, $|\text{{\it \underline{Order-Types}}}(J(p))|=\mathcal{O}(p^2)$.

For arbitrary $k$ one would start with the configuration that has $i_1,\dots, i_{k+1}$ at the vertices of the regular $(k+1)$-gon, then placing a pair of points $b_j^1,b_j^2$ on the arc connecting $i_j$ and $i_{j+1}$. Finally, there is $k$ points $x_1,\dots,x_k$ of the inner layer, placed close enough to the center of each segment $[b_j^1,b_j^2]$.
One makes sure that $x_j$ is in triangle $T(t,b_j^1,b_j^2)$, where $t$ ranges over all points of the outside layer, other than $b_j^1,b_j^2$. In particular, all elements $i_1,\dots, i_{k+1}$ are equivalent. As for examples above, place $p>k$ points into locations of $i_1,\dots, i_{k+1}$, splitting $p$ into subsets with $p_1,\dots,p_k,(p-(p_1+\dots+p_k))$ points, correspondingly, for some $p_1+\dots+p_k \leq p$. This produces a point realization for a convex geometry $J(p)$ from series $\mathcal{L}_k$, whose size grows as $\mathcal{O}(p)$, while the number of non-equivalent order types is $\mathcal{O}(p^k)$. 
\end{proof}

\begin{cor}\label{plenty} The growth of $|\text{{\it \underline{Order-Types}}}(J(p))|$ of two-layered convex $4$-geometries $J(p)$ of size $\mathcal{O}(p)$ cannot be $p$-polynomially bounded.
\end{cor}

\section{Quasi order types}

Suppose $t = t_f$ is an order type on a set $J$, and $x,y \in J$ are such that either
\vskip .4cm

(3) \quad  $(\forall z \in J - \{x,y\} ) \ \  t (x,y,z) = 1$\quad or\quad $(\forall z \in J- \{x,y\}) \ \ t (x,y,z) = -1$.

\vskip .4cm

Then obviously $f(x), f(y)$ must be two adjacent points of the outermost layer of the point configuration $f(J) \subseteq \mathbb{R}^2$. The converse holds as well. This motivates the following concept. For any function $t: J[3] \rightarrow \{-1,1\}$, call $\{x,y\} \subseteq J$ a {\it quasi-edge} of $t$ if $(3)$  takes place. The set $J$ being finite we may recursively define a {\it quasi order type} as any function $t : J [3] \rightarrow \{-1, 1\}$ such that
\begin{enumerate}
 \item [(i)] the graph $G$, whose edge set $E(G)$ is the set of all quasi-edges of $t$, is a cycle (called a {\it quasi layer}), and
\item[(ii)] the set $J': = J- \cup E(G)$ is either empty, or the restriction of $t$ to $J'[3]$ is a quasi order type.
\end{enumerate}
It is clear that each order type is a quasi order type in such a way that layers and quasi layers coincide. Furthermore, it takes time $O(n^3) \ (n = |J|)$ to check whether or not a function $t: \ J[3] \rightarrow \{-1,1\}$ is a quasi order type.

The \emph{depth} of a quasi order type is the number $q(t)$ of its quasi layers.

A {\it quasi rooted triangle} of a quasi order type $t: \ J[3] \rightarrow \{-1,1\}$ is a pair $(\{a, b, c\},d)$ satisfying 
$$t(a,b,c) = t(d, b, c) = t(a,d,c) = t(a,b,d)$$
It takes time $O(n^4)$ to compute the set $\mathcal{Q} \mathcal{R} \mathcal{T}(t)$ of all quasi rooted triangles of $t$. Dito it costs $O(n^4)$ to decide whether $\mathcal{Q}\mathcal{R} \mathcal{T}(t)$ satisfies the  Dietrich axiom and hence yields a convex $4$-geometry [Prop.2.1]. A quasi order type $t$ satisfying the Dietrich axiom is {\it simple} if $\mathcal{Q}\mathcal{R}\mathcal{T}(t)$ yields a simple convex $4$-geometry in the sense of Definition \ref{simple}.

\begin{rmk}
Observe that it can be tested in polynomial time whether or not a function $t : \ J[3] \rightarrow \{-1,1\}$ is a simple quasi order type.
\end{rmk}

%\vskip 1cm

\section{Complexity of the modified Edelman-Jamison problem} 

We tempt to link the Edelman-Jamison Problem \ref{EJP} to
\begin{pbm}
{\bf The Order Type Problem}
Given any function $t: J[3] \rightarrow \{1, -1\}$, recognize whether it is an order type and, if it is, find some realizing point configuration.
\end{pbm}

It is known that the Order Type Problem is NP-hard; that follows from the famous Mn\"{e}v's Universality Theorem \cite{Mn}. 

In this section we consider the modified Edelman-Jamison Problem and we will show that it is polynomially equivalent to the Order-Type Problem.

\begin{pbm}\label{MEJP}
Suppose we are given a convex $4$-geometry $(J,\mathcal{R}\mathcal{T})$ and, in addition, some fixed circular order of the outermost layer $L_0=\{a_1,\dots,a_n\}$.
{\bf The modified Edelman-Jamison Problem} asks whether this geometry can be realized by a point configuration in the Euclidean plane {\bf with this given clock-wise circular order of the outermost layer}.
\end{pbm}
We will say that a function $t^*: J[3] \rightarrow \{1, -1\}$ supports the clock-wise circular order $L_0=\{a_1,\dots,a_n\}$, if $t^*(a_i,a_j,a_k)=-1$ for all 
$i<j<k$ (modulo $n$). 

\begin{thm}\label{main}
Given a finite convex $4$-geometry $(J,\mathcal{R}\mathcal{T})$ and a circular order $L_0
=\{a_1,\dots,a_n\}$ of its outermost layer, one can decide in polynomial time that either this geometry is not realizable with such circular order, or define a unique 
%up to the opposite,
function $t^*: J[3] \rightarrow \{1, -1\}$ associated with the geometry that supports the clock-wise circular order $L_0$. If $(J, \mathcal{R}\mathcal{T})$ happens to be realizable, then every point realization of $(J,\mathcal{R}\mathcal{T})$ with the given clock-wise order of $L_0$ will be, as an order type, weakly equivalent to $t^*$.
\end{thm}

In other words, if the geometry is realizable then knowing the clock-wise circular order of its outermost layer defines uniquely the order type of its realization.

In order to prove Theorem \ref{main}, we need to introduce the following definitions.

\begin{df}\label{carousel}
A circular ordering $\{a_1,\dots,a_n\}$ of the outermost layer $L_0$ of a convex $4$-geometry $(J, \mathcal{R}\mathcal{T})$ satisfies  \emph{the carousel rule} if, for any element $x \in J$, and any element $y \in J\backslash L_0$, 
there exists exactly one $i\leq n$ such that
\begin{equation}
(\{x, a_i, a_{i+1} \}, y) \in \mathcal{R}\mathcal{T} \quad (\mbox{modulo} \ n)
\tag{CR}
\end{equation}
\end{df}

\vspace{0.5 cm}

\begin{figure}[ht]
\begin{center}
\includegraphics[scale=.45]{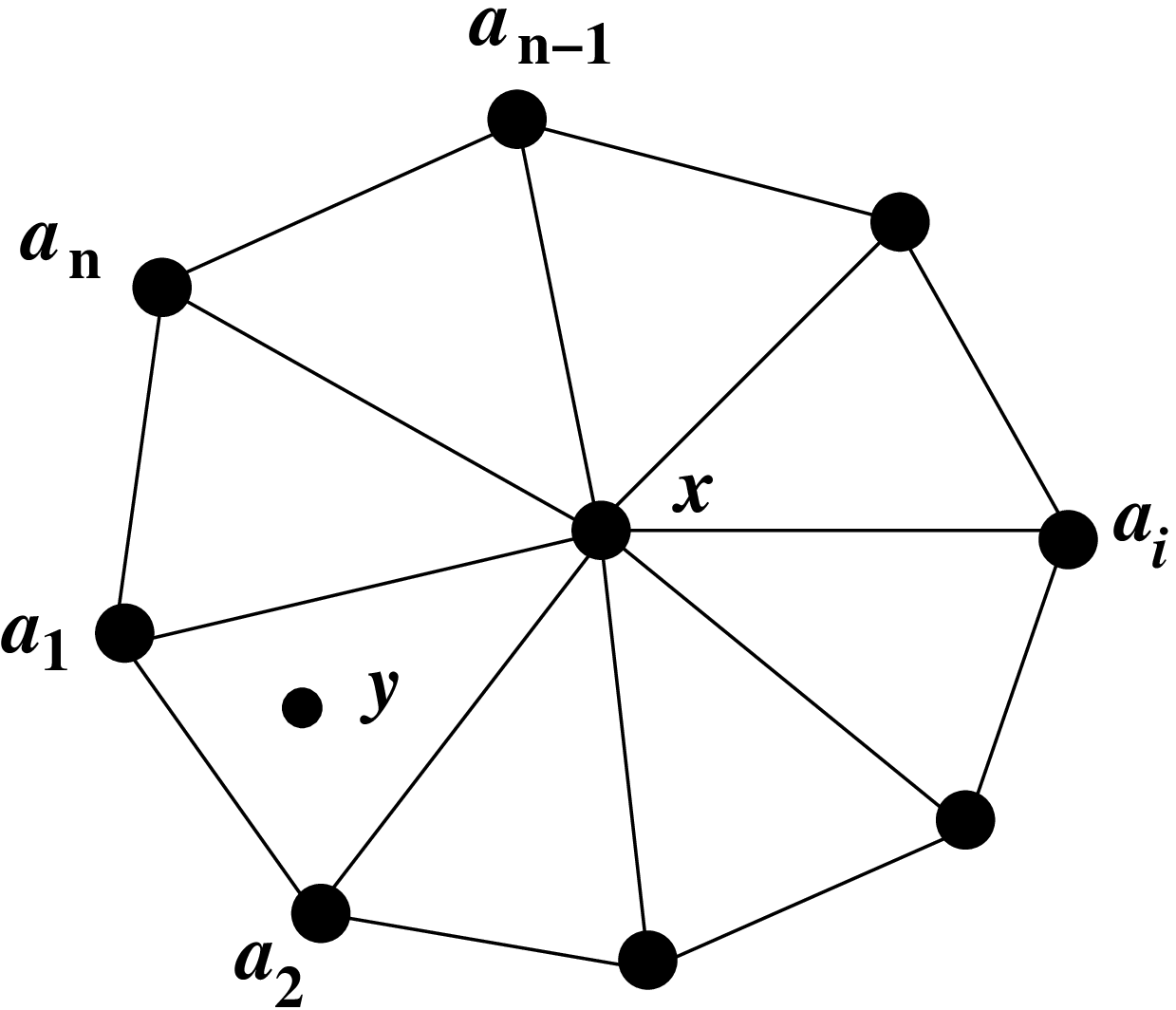}
\caption{}
\label{Fig. 2}
\end{center}
\end{figure}

Figure \ref{Fig. 2} illustrates the carousel rule in the realizable case. Point $y$ belongs to the polytope with vertices $a_1, \dots, a_n$, in that circular order, but it is not a vertex of this polytope. 
For any other point $x$ of that polytope (including the case
when $x$ is one of $a_1,\dots, a_n$) we consider the splitting
of the polytope into triangles $\{x,a_i,a_{i+1}\}$. The carousel rule is the statement that
$y$ belongs to only one of those triangles.
\newpage

\begin{df}\label{3-carousel}
We will say that the convex $4$-geometry $(J,\mathcal{R}\mathcal{T})$ satisfies the $3$-carousel rule,
if for any five distinct points $a, b, c, x, y \in J$ with $x,y \in \overline{\{a, b, c\}}$ exactly one of the following alternatives takes place:
$$x \in \overline{\{a, b, y\}}, \quad x \in \overline{\{a, c, y\}} \quad x \in \overline{\{b,c,y\}}$$
Of course, the first statement, say, amounts to $(\{a, b, y \}, x) \in \mathcal{R}\mathcal{T}$, but this notation would be a bit clumsy.
\end{df}

\begin{lm}\label{necessary} Every realizable $4$-geometry $(J,\mathcal{R}\mathcal{T})$ satisfies $3$-carousel rule. Besides, if $\{a_1, \dots, a_n\}$ is a circular ordering of the outermost layer in some point realization, then this ordering satisfies the carousel rule.
\end{lm}
The proof is evident.

\begin{lm}\label{check_carousel}
Given a finite convex $4$-geometry $(J,\mathcal{R}\mathcal{T})$ and a circular order $L_0
=\{a_1,\dots,a_n\}$ of its outermost layer, it can be checked in polynomial time whether
this ordering satisfies the carousel rule, and this convex $4$-geometry satisfies $3$-carousel rule.
\end{lm}

\begin{lm}\label{easy}
In a convex $4$-geometry that satisfies $3$-carousel rule the following holds:
if $b,x \in T(a,a_j,a_{j+1})$ and $x \in T(b,a_j,a_{j+1})$ then one and only one statements holds: $b \in T(x,a,a_j)$ or $b \in T(x,a,a_{j+1})$.
\end{lm}
\begin{proof} By the $3$-carousel roule the only other possiblity for $b$ a priori is $b \in \overline{\{a_j, a_{j+1}, x\}}$. However, because of $x \in \overline{\{a_j, a_{j+1}, b\}}$ and the anti-exchange property, this can't happen.
\end{proof}

\begin{proof} of {\it Theorem \ref{main}}.
Suppose we are given a finite convex $4$-geometry $(J,\mathcal{R}\mathcal{T})$ and a circular order $L_0
=\{a_1,\dots,a_n\}$ of its outermost layer. Due to Lemma \ref{check_carousel}, in polynomial time one can check whether the carousel rule is satisfied for the given ordering of $L_0$, and whether the $3$-carousel rule holds. If either fails, then, according to Lemma \ref{necessary}, the geometry is not realizable with the given ordering of $L_0$.

If they both hold then we are going to define the unique function $t^*:J[3] \rightarrow \{1, -1\}$ associated with the given geometry and supporting ''clock-wise'' ordering of $L_0$.

In order to define $t^*:J[3] \rightarrow \{1, -1\}$ we will show that every ordered pair $(a,b) \in J^2$ triggers a unique splitting of $J_1=J \backslash \{a,b\}$ into two subsets $K$ and $K'$.  We will be guided by the requirement that, in case the geometry happens to be realizable by some point configuration where the points $a_1,\dots,a_n$ of the outermost layer follow clock-wise, the proposed splitting will represent the splitting of $J_1$ into these two subsets: those points that lie in the ''left'' semi-plane with respect to (the suitably directed) $line(a,b)$, and those that lie in the ``right'' semi-pane. Thus, we would define $t^*(a,b,x)=1$ for every $x \in K$ and $t^*(a,b,x)=-1$ for every $x \in K'$.

In fact, referring to the definition of the order type \cite{GoPo}, the knowledge of how a set of points in $\mathbb{R}$ is split into two subsets $K$ and $K'$ (''left'' and ''right'') by any line through two points, defines the order type of the given configuration up to weak equivalence.

If $a,b$ are points in $L_0$, then assuming, say, $a=a_i$, $b=a_j$ with $i<j$, we define $K=J_1 \cap \overline{\{a_i,a_{i+1},\dots,a_j\}}$ and $K'=J_1 \cap \overline{\{a_j,a_{j+1},\dots,a_i\}}$ (modulo $n$). It follows from the carousel rule for this ordering of $L_0$ that every point in $J_1$ will be exactly in one of $K$ or $K'$.

If $a \in L_0$ and $b \not \in L_0$, then $a$ is, say, $a_1$, and, according to the carousel rule, $b \in T := \overline{\{a,a_j,a_{j+1}\}}$, for uniquely defined $j >1$. Again, due to carousel rule, every point of $J_1$ will be exactly in one of three sets: $A_1=\overline{\{a_1,\dots,a_j\}}$, $T$, or $A_2=\overline{\{a_j, \dots, a_n,a_1\}}$. Besides, every point of $J_1$ that gets into $T$, will be exactly in one of three sets:
$$T_1 : = \overline{\{a, a_j, b\}}, \quad T_2 := \overline{\{a, a_{j+1}, b\}}, \quad T_3 := \overline{\{b, a_j, a_{j+1} \}}$$
Thus, every point of $J_1$ will be in one and only one of these sets:
$$B_1: = A_1 \cup T_1, \quad B_2 : = A_2 \cup T_2, \quad T_3$$
By Lemma 6.8, for each $x \in T_3$ exactly one of these statements is true: $b \in \overline{\{x, a, a_j \}}$ or $b \in \overline{\{x, a, a_{j+1} \}}$. Hence it is clear that we need to define
$$\begin{array}{lll} K & : = & (B_1 \cup \{x \in T_3 | \ b \in \overline{\{x, a, a_{j+1} \} } \} ) \setminus \{a, b\} \\
\\
K' & := & (B_2 \cup \{x \in T_3 | \ b \in \overline{\{x, a, a_j \}} \} \ ) \setminus \{a, b\}
\end{array}$$
The last case, when both $a,b$ are not in $L_0$, is similar to the previous. One finds, due to the carousel rule, the unique $i<j$ such that, say, 
$a \in \overline{ \{b, a_i, a_{i+1} \}}$, and $b \in  \overline{\{ a,a_j,a_{j+1}\}}$. Then every point of $J_1\setminus L_0$ will be in one and only one of the following four sets (Fig.8 visualizes the realizable case):
$$T: = \overline{ \{ a, a_i, a_{i+1} \}}, \ A_1 := \overline{\{a, a_{i+1}, \cdots, a_j \}}, \ T' : = \overline{\{a, a_j, a_{j+1} \}}, \ A_2: = \overline{\{a, a_{j+1}, \cdots, a_i \}}$$

\newpage
\begin{figure}[!h]
\begin{center}
\includegraphics[scale=.45]{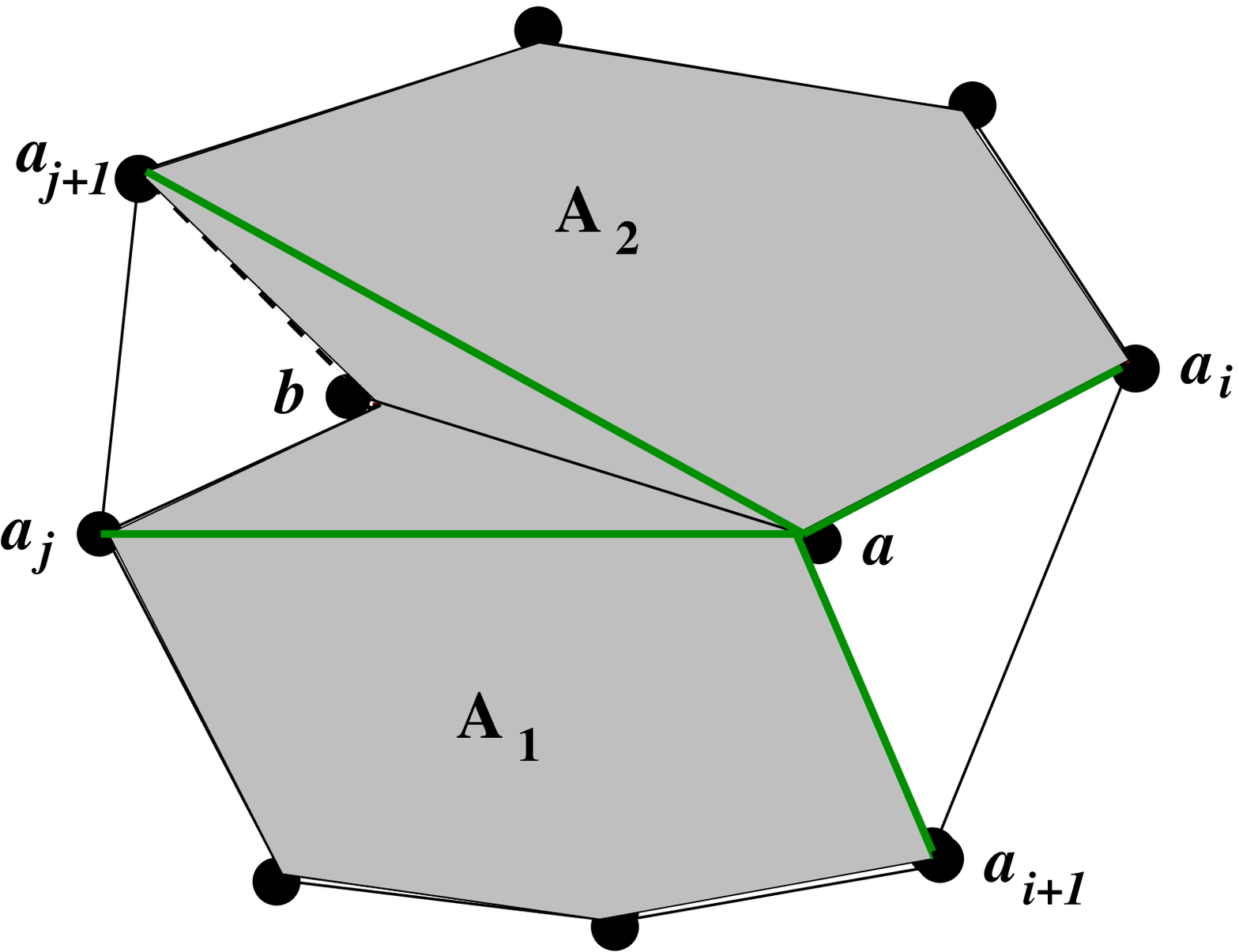}
\caption{}
\end{center}
\end{figure}

Due to the $3$-carousel rule every point of $T' \setminus \{a, a_j, a_{j+1}, b\}$ is in one and only one of the sets
$$T'_1 : = \overline{\{b, a, a_j\}}, \ T'_2 := \overline{\{b, a_j, a_{j+1} \}}, \ T'_3 := \overline{\{b, a_{j+1}, a\}}$$
Therefore, if we put
$$\begin{array}{lll} K & := & (A_1 \cup T'_1 \cup \{x \in T'_2 | \ b\in \overline{\{x,a,a_{j+1}\}} \} \cup \{ x \in T| \ a \in \overline{\{x, b, a_i\}} \} ) \setminus \{a, b\} \\
\\
K'  & : = & (A_2 \cup T'_3 \cup \{x \in T'_2 | \ b \in \overline{\{x, a, a_j \}} \} \cup \{x \in T | \ a \in \overline{ \{x, b, a_{i+1} \}}) \setminus \{a, b\} \end{array}$$
then $K, K'$ is a bipartition of $J_1$ which, should $(J, \mathcal{R}\mathcal{T})$ be realizable, is induced by the line through $a$ and $b$.

Assume now that given convex $4$-geometry $(J,\mathcal{R}\mathcal{T})$ is realizable with the clock-wise order $L_0$ of its outermost layer. Then the procedure described above corresponds to splitting points realizing $J$ into two subsets, for any points $a,b \in J$: to those that belong to the ''left'' semi-plane with respect to directed $line (a,b)$, and to those that belong to ''right'' semi-plane. Thus, produced function $t^*:J[3] \rightarrow \{1, -1\}$ represents an order-type of the configuration. Besides, it supports the clock-wise order of $L_0$. It proves that all possible point configurations of $(J,\mathcal{R}\mathcal{T})$ with the given clock-wise order of $L_0$ will be equivalent as order-types.
\end{proof}

\begin{cor}\label{equiv}
The modified Edelman-Jamison Problem is equivalent to the Order-Type Problem. In particular, the modified Edelman-Jamison Problem is NP-hard.
\end{cor}
\begin{proof}
Indeed, given an instance of modified Edelman-Jamison Problem, in polynomial time one can check whether $(J, \mathcal{R}\mathcal{T})$ it satisfies carousel rule for the given order of the outermost layer, and whether $(J, \mathcal{R}\mathcal{T})$ satisfies $3$-carousel rule. If not, $(J, \mathcal{R}\mathcal{T})$ is not realizable. Otherwise, we obtain a uniquely defined function $t^*:J[3] \rightarrow \{1, -1\}$ that supports the given clock-wise order of the outermost layer, thus, obtaining the instance of the Order Type Problem. If the latter is solved positively, the same point configuration will provide the solution to modified Edelman-Jamison Problem. If it is solved negatively, the modified Edelman-Jamison Problem is refuted, too.

Vice versa, given an instance of the Order Type Problem, in polynomial time, one can either refute, or assert that the given function $t:J[3] \rightarrow \{1, -1\}$ is a quasi order type. In particular, the convex $4$-geometry will be defined together with the clock-wise ordering of the outermost layer, thus, we will get the instance of the modified Edelman-Jamison Problem. According to Theorem \ref{main}, if this convex geometry with the given clock-wise ordering is realized, then such a realization will represent a unique order-type $t^*$ that supports the given clock-wise order of the outermost layer. Note that $t$ and $t^*$ will agree on any triple of elements from $L_0$. It will take polynomial time to check whether $t$ and $t^*$ agree on all triples, and if they do, the Order Type Problem is solved positively, otherwise - not.

If the geometry cannot be realized with the given ordering of the outermost layer, then the quasi order type is not an order type. 
\end{proof}

\begin{cor}
The following problems are equivalent:
\begin{enumerate}
\item [(a)] There is a polynomial time algorithm which decides whether a simple quasi order type $t:  \ J[3] \rightarrow \{-1,1\}$ is an order type. 
\item[(b)] There is a polynomial time algorithm that decides whether a given simple convex geometry is realizable.
\end{enumerate}
\end{cor}
\begin{proof}
This follows from Corollary \ref{equiv} due to Theorem \ref{ordering}. Indeed, simple convex $4$-geometry has a unique circular ordering of the outermost layer. Thus, the modified Edelman-Jamison Problem for such geometry is equivalent to Edelman-Jamison Problem.

\end{proof}

{\bf Acknowledgments.} The results of this paper were presented on the geometry seminar at Courant Institute of Mathematical Sciences, in New York, in spring of 2006. The results were also presented at Colloquium of the Mathematics Department of Iowa State University, in February of 2007. We are grateful to Prof.R.Pollack, who organizes the geometry seminar, and Prof.J.D.H. Smith and Prof. A.Romanowska, faculty at Iowa State, for their interest in our results, and for arrangements of the seminar visits of the first author. We appreciate the help of Fedor Adarichev who translated part of the data from \cite{AAK} to a printable image of order types that we used in our study. We were helped by Vyacheslav Adarichev and PhD student Yves Semegni in preparing pictures for the paper.

\end{document}